\documentclass[11pt,reqno]{amsart}
\usepackage{ amsmath,  amsfonts, amssymb,  amsthm,   amscd, gensymb,  graphicx, comment,  etoolbox, url,
booktabs, stackrel, mathtools,enumitem, mathdots,  microtype, lmodern,    mathrsfs, graphicx, tikz,  longtable,tabularx, float, tikz, pst-node, tikz-cd, multirow, tabularx, amscd,  bm, array, makecell, diagbox, booktabs,ragged2e, caption, subcaption }
\usepackage{makecell,slashbox}
\usepackage{xcolor}
\usepackage[all]{xy}
\usepackage[pagebackref=true, colorlinks=true, linkcolor=blue, citecolor=blue, urlcolor=blue, breaklinks=true]{hyperref}
\usepackage[utf8]{inputenc}
\usepackage{microtype, fullpage, wrapfig,textcomp,mathrsfs,csquotes,fbb}
\usepackage[colorlinks=true, linkcolor=blue, citecolor=blue, urlcolor=blue, breaklinks=true]{hyperref}
 
\usepackage[capitalise]{cleveref}
\usepackage{todonotes}
\usetikzlibrary{positioning}
\usetikzlibrary{shapes,arrows.meta,calc}
\usetikzlibrary{arrows}

\newtheorem{theorem}{Theorem}[section]

\newtheorem{lemma}[theorem]{Lemma}
\newtheorem{point}[theorem]{}

\newtheorem{proposition}[theorem]{Proposition}
 
\newtheorem*{theorem*}{Theorem}
\newtheorem*{conjecture*}{Conjecture}
\theoremstyle{definition}
\newtheorem{definition}[theorem]{Definition}
\newtheorem{example}[theorem]{Example}
\newtheorem{remark}[theorem]{Remark}

\newtheorem{question}[theorem]{Question}
\newtheorem{setup}[theorem]{Setup}

\newtheorem*{notation*}{Notation}
 
\newcommand{\p}{\mathfrak{p}}

\newcommand{\m}{\mathfrak{m}}

\newcommand{\Z}{\mathbb{Z}}
\newcommand{\N}{\mathbb{N}}

\newcommand{\Ass}{\operatorname{Ass}}

\newcommand{\Hom}{\operatorname{Hom}}

\newcommand{\Ext}{\operatorname{Ext}}
\newcommand{\Tor}{\operatorname{Tor}}

\makeatletter
\def\@tocline#1#2#3#4#5#6#7{\relax
  \ifnum #1>\c@tocdepth % then omit
  \else
    \par \addpenalty\@secpenalty\addvspace{#2}%
    \begingroup \hyphenpenalty\@M
    \@ifempty{#4}{%
      \@tempdima\csname r@tocindent\number#1\endcsname\relax
    }{%
      \@tempdima#4\relax
    }%
    \parindent\z@ \leftskip#3\relax \advance\leftskip\@tempdima\relax
    \rightskip\@pnumwidth plus4em \parfillskip-\@pnumwidth
    #5\leavevmode\hskip-\@tempdima
      \ifcase #1
       \or\or \hskip 1em \or \hskip 2em \else \hskip 3em \fi%
      #6\nobreak\relax
    \hfill\hbox to\@pnumwidth{\@tocpagenum{#7}}\par% <---- \dotfill -> \hfill
    \nobreak
    \endgroup
  \fi}
\makeatother
\newcolumntype{x}[1]{>{\centering\arraybackslash}p{#1}}

\begin{document}
\title[]{A study on v-numbers of $\Ext^i_R(L,M/I^nN)$ and $\Tor_i^R(L,M/I^nN)$}

\author[Sayed Sadiqul Islam]{Sayed Sadiqul Islam}
\address{Indian Institute of Technology Bombay, Mumbai, India}
\email{22d0786@iitb.ac.in}

\author{Lisa Mandal}
\address{Department of Mathematics and Statistics, Indian Institute of Science Education and Research Kolkata, Mohanpur, India}
\email{lishamandal22@gmail.com}

\author{Samarendra Sahoo}
\address{Indian Institute of Technology Dharwad, Dharwad 580011, Karnataka, India}
\email{samarendra.s.math@gmail.com, ra.samarendra@iitdh.ac.in}

\thanks{}
 
\begin{abstract} 
Let $R$ be a Noetherian $\N$-graded ring. Assume that $L$ and  $M$ are finitely generated $\Z$-graded $R$-module, and $I$ a homogeneous ideal of $R$. Let $N$ be a graded submodule of $M$. Fix $k\in \N$. Denote $H:=\Ext^k_R(L,M)$ (or $\Tor_k^R(L,M)$) and $H_n:=\Ext^k_R(L, M/I^nN)$ (or $\Tor_k^R(L,M/I^nN)$). When $I$ is generated by elements of positive degree, we show that the local $v$-number of $H_n$ either coincides with that of colon submodule $(0:_HI)$ or is a polynomial in $n$ of degree one whose leading coefficient is one of the degrees of the generators of $I$. This extends the main result of Ghosh-Pramanik [J. Algebra 671 (2025), 61–74], where asymptotic linearity was shown under the assumption that $(0:_HI)=0$.
\end{abstract}
\keywords{Asymptotic associated primes; Graded rings and modules; Ext, Tor}
\subjclass[2020]{Primary 13D07, 13A02,  13A15}
\maketitle
\section{Introduction}
The concept of the Vasconcelos invariant, associated with a homogeneous ideal in a polynomial ring over a field, was proposed in \cite{CST+20}  to investigate the limiting behavior of the minimum distance of projective Reed-Muller type codes. This invariant, currently known as the $v$-number, was generalized to the setting of graded modules in \cite{FG25}. The asymptotic behavior of $v$-numbers for graded modules has attracted significant attention in recent years (see, for instance, \cite{BMS26, FS25, VS25, FM25, VS24,  Con24, Fic24, FS26}). In addition to these developments, our investigation is motivated by several other prominent works \cite{FG25, GP25, GNP26} in this direction.
\begin{setup}\label{main setup}
     Let $R$ be a Noetherian $\N$-graded ring. Assume that $L$ and  $M$ are finitely generated $\Z$-graded $R$-module, and $I$ a homogeneous ideal of $R$. Let $N$ be a graded submodule of $M$. Suppose $I$ is generated by homogeneous elements $y_1,\ldots,y_c$ of degree $d_1\leq\ldots\leq d_c$ respectively.
\end{setup}
 Recall, for $\p \in \Ass_R(M)$, the  local $v$-number of $M$ at $\p$ is defined as $v_{\p}(M)= \inf\{n\, |\, \p = (0 :_R x)\text{ for some } x \in M_n\}.$
The global $v$-number of $M$ is defined by $$v(M)=\min\left\{v_{\p}(M):\p\in\Ass_R(M)\right\}.$$
By convention, $v(0)=\infty$. 
In \cite{Bro79}, Brodmann proved that the set $\Ass_R(M/I^nM)$ stabilize eventually. Let $\mathcal A_M(I)$ be the stable set. Motivated by this result of Brodmann, Conca showed that when $R$ is a domain and $\p\in \mathcal A_R(I)$, $v_{\p}(R/I^n)=a_\p n+b_\p$  for all $n\gg 0$, where $a_p$ is equal to one of the degrees of the generators of $I$, and $b_\p\in \Z$ \cite[Theorem 1.1]{Con24}. When $R$ is a polynomial ring over a field, the same result was proved by Ficarra-Sgroi in \cite{FS26}. In \cite[Theorem 2.14]{FG25}, Fiorindo-Ghosh showed that when  $\p\in \mathcal A_M(I)$, then for $n\gg 0$, $v_{\p}(M/I^nM)$  is linear functions of $n$ whose leading coefficients are from the degree of the generators of $I$ provided $(0:_MI)=0$. In a recent paper \cite{GNP26}, the authors study the case when $(0:_MI)\neq 0$. When $I$ is generated by elements of positive degree, the authors  in \cite{GNP26} show that the local $v$-number of $M/I^nM$ either coincides with that of colon submodule $(0:_MI)$ or is a polynomial in $n$ of degree one whose leading coefficient is one of the degrees of the generators of $I$.

Now let $M$, $N$ and $L$ be finitely generated graded $R$-modules with $N\subseteq M$. In \cite[Corollary 3.5]{KW04}, Katz-West proved that, for fixed $k\in \N$, the sets $\Ass_R(\Ext^k_R(L,M/I^nN))$ and  $\Ass_R(\Tor_k^R(L,M/I^nN))$ are independent of $n$ for $n\gg 0$.  In \cite{GP25}, Ghosh-Pramanik proved that, for fixed $k\in \N$, the local $v$-number of $\Ext^k_R(L, M/I^nN)$ is eventually a linear function of $n$ for $n\gg0$ with leading coefficient being the degree of a generators of $I$ whenever $(0:_{\Ext^k_R(L,M)}I)=0.$ They also showed that an analogous statement holds for $\Tor_k^R(L,M/I^nN)$  whenever $(0:_{\Tor_k^R(L,M)}I)=0.$ Keeping these results in mind, the following questions arise naturally.
\begin{question}
Let $R$ be a Noetherian $\N$-graded ring. Assume that $L$ and  $M$ are finitely generated $\Z$-graded $R$-module, and $I$ a homogeneous ideal of $R$. Let $N$ be a graded submodule of $M$. Fix $k\in \N$. Denote $H:=\Ext^k_R(L,M)$ (or $\Tor_k^R(L,M)$) and $H_n:=\Ext^k_R(L, M/I^nN)$ (or $\Tor_k^R(L,M/I^nN)$). Let $(0:_HI)\neq 0$.
\begin{enumerate}
    \item Let $\p\in \Ass_R(H_n)$ for $n\gg 0$. What is the asymptotic behaviour of $v_\p(H_n)$ as a function of $n?$
    \item What is the asymptotic behaviour of $v(H_n)?$
\end{enumerate}
\end{question}
We answer the above question via the following theorem which is the main result of this paper.
\begin{theorem}(Theorem \ref{proof of main result} below)\label{main result}
   Let $R$ be a Noetherian $\N$-graded ring. Assume that $L$ and  $M$ are finitely generated $\Z$-graded $R$-modules, and $I$ a homogeneous ideal of $R$,  generated by homogeneous elements $y_1,\ldots,y_c$ of degree $d_1\leq\ldots\leq d_c$ respectively. Let $N$ be a graded submodule of $M$. Fix $k\in \N$. Denote $H:=\Ext^k_R(L,M)$ (or $\Tor_k^R(L,M)$) and $H_n:=\Ext^k_R(L, M/I^nN)$ (or $\Tor_k^R(L,M/I^nN)$). Then the following statements hold:
   \begin{enumerate}[label=(\arabic*), font=\upshape]
       \item  Let $\p\in \left(\Ass_R H_n\cap V(I)\right)\setminus \Ass_R(0:_HI)$ for $n\gg 0$. Then there exists $a_\p\in\{d_1,\ldots,d_c\}$ and $b_\p\in \Z$ such that $v_{\p}(H_n)=a_\p n+b_\p$ for $n\gg 0$.
       \item  Let $\p\in \Ass_R(H_n)\cap \Ass_R(0:_HI)$ for $n\gg 0$. If $d_i\geq 1$ for $1\leq i\leq c$, then $v_{\p}(H_n)=v_{\p}(0:_HI)$ for $n\gg 0$.
       \item Let $(0:_HI)\neq 0$. Then, $v(H_n)\leq v(0:_HI)$ for $n\gg 0$. Further, if $d_i\geq 1$ for $1\leq i\leq c$, and $\Ass_R(H_n)\subseteq V(I))$, then $v(H_n)= v(0:_HI)$ for $n\gg 0$. 
   \end{enumerate}
\end{theorem}
\textbf{Organization of the article:} 
This paper is structured into five sections. Section 2 introduces the required preliminaries. In Section 3, we establish the formula for $\Ass_R \left(U+I^nV/I^nW\right)$, while Section 4 presents the proof of our main result (Theorem \ref{main result}). Section 5 concludes with examples demonstrating our main results.

\section{Preliminaries}
Throughout this article, all rings are commutative Noetherian rings with unity, and all modules are finitely generated graded $R$-modules. In this section, we discuss some preliminary results concerning associated primes and local $v$-numbers of graded modules that will be useful in our proofs.
\begin{definition}
    The \emph{$I$-torsion submodule} of  an $R$-module $M$ is defined by
$\Gamma_I(M)
=\{m\in M : I^n m=0 \text{ for some } n\geq 1\}
=\bigcup_{n\geq 1}(0:_M I^n).$
\end{definition}

\begin{lemma}\cite[Lemma 2.5]{GNP26}\label{ass prime of gammaI}
For any $R$-module $M$ and an ideal $I$ of $R$ we have     $\Ass_R( M)\cap V(I)=\Ass_R(\Gamma_I(M))=\Ass_R((0:_RI^n))$ for all $n\geq 1$.
\end{lemma}
We frequently rely on the next two propositions throughout our proofs. The first one shows the behavior of $v$-numbers in an exact sequence.
\begin{proposition}\cite[Proposition 3.1]{GNP26}\label{v number in exact seq}
    Let $0\to M_1\to M_2 \to M_3\to 0$ be a short exact sequence of finitely generated $\Z$-graded modules over a Noetherian $\N$-graded ring $R$. Then the following statements are true.
    \begin{enumerate}[label=(\arabic*), font=\upshape]
        \item If $\p\in \Ass( M_1)$, then $v_{\p}(M_2)\leq v_{\p}(M_1).$ Consequently, $v(M_2)\leq v(M_1).$ 
        \item If $\p\in \Ass_R (M_2)\setminus\Ass_R (M_1)$, then $v_{\p}(M_3)\leq v_{\p}(M_2).$
    \end{enumerate}
\end{proposition}
The following proposition relates the local $v$-numbers of a module, its $I$-torsion submodule, and the associated colon ideals $(0:_MI^n)$, which will be crucial for our subsequent proofs.
\begin{proposition}\cite[Proposition 3.2]{GNP26}\label{v number of M and colon same}
    Let $R$ be a Noetherian $\N$-graded ring, $M$ a finitely generated $\Z$-graded $R$-module, an $I$ a homogeneous ideal of $R$.
    \begin{enumerate}[label=(\arabic*), font=\upshape]
        \item If $\p\in \Ass_R((0:_MI))$, then $v_{\p}(M)=v_{\p}(\Gamma_I(M))=v_{\p}(0:_MI^n)$ for each $n\geq 1$.
        \item $v(\Gamma_I(M))=v(0:_MI^n)$ for each $n\geq 1$.
        \item If $\Ass_R(M)\subseteq V(I)$, then $v(M)=v(\Gamma_I(M))$.
    \end{enumerate}
\end{proposition}
The following result is an elementary observation concerning the $v$-numbers.

\begin{proposition}\label{easy}
Let $M$ and $N$ be finitely generated
$\mathbb{Z}$-graded $R$-modules. Then, 
\begin{enumerate}
    \item For every
$\p \in \Ass_R (M)\cap \Ass_R (N)$,
we have $v_{\p}(M\oplus N)=\min\{v_{\p}(M),v_{\p}(N)\}.$
\item Fix $r\in \Z$. We have $v(M(r))=v(M)-r$ and $v_{\p}(M(r))=v_{\p}(M)-r$ for all $\p \in \Ass_R(M).$
\end{enumerate}
\end{proposition}

The following well-known Artin-Rees Lemma will be useful in our later proofs.
\begin{lemma}\cite[Theorem 17.1.6]{HS06}\label{Artin-Ress}
    Let $R$ be a Noetherian ring, $I$ an ideal of $R$, and $M$ be a finitely generated $R$-module. Let $N$ be a submodule of $M$. Then there exists an integer $c\geq 0$ such that 
    $$I^mM\cap N =I^{n-c}(I^cM\cap N),\qquad  \text{for all $n\geq c.$}$$
\end{lemma}

We now recall the following two lemmas that says $\Ext^k_R(L,M/I^nN)$ and $\Tor_k^R(L,M/I^nN)$ can be expressed as a quotient of graded modules involving powers of $I$.
\begin{lemma}\cite[Lemma 3.3]{GP25}\label{ext lemma}
    With Setup \ref{main setup}, fix $k\in \N$. Then there exists $q\in \N$ such that
    $$\Ext^k_R(L,M/I^nN)\cong (U+I^{n-q}V)/I^{n-q}W \ \text{for all}\ n\geq q$$
    and for some graded submodules $U,V$ and $W$ of a finitely generated $R$-module $E$ satisfying $W\subseteq V$ and $U\cong \Ext^k_R(L,M)$.
\end{lemma}
An analogous result holds for Tor-modules, as stated in the following lemma.
\begin{lemma}\cite[Lemma 3.4]{GP25}\label{tor lemma}
    With Setup \ref{main setup}, fix $k\in \N$. Then there exists $q\in \N$ such that
    $$\Tor_k^R(L,M/I^nN)\cong (U+I^{n-q}V)/I^{n-q}W \ \text{for all}\ n\geq q$$
    and for some graded submodules $U,V$ and $W$ of a finitely generated $R$-module $E$ satisfying $W\subseteq V$ and $U\cong \Tor_k^R(L,M)$.
\end{lemma}

Throughout this article, we assume the following hypothesis

\begin{point}\label{hypo}
    \normalfont
\textbf{Hypothesis:}  Let $R$ be a Noetherian $\N$-graded ring and $I$ be a homogeneous ideal of $R$ generated by homogeneous elements $y_1,\ldots, y_c$ of degree $d_1\leq \ldots \leq d_c$ respectively. Let $U,V,W$ be submodules of a finitely generated $\Z$-graded $R$-module $E$ such that $W\subseteq V$. 

For each integer $n\geq 1$, put
$M_n=U+I^nV, D_n=
M_n\cap\left((0:_E I)+I^{n-1}W\right), M_n'=M_n/I^nW, D_n'=D_n/I^nW,$ and $T_n=(M_n\cap I^{n-1}W)/I^nW$. Note that $I^nW\subseteq M_n$ as $W\subseteq V$. Also since $I^nW\subseteq D_n$, 
$D_n'$ is well-defined.
\end{point}

Note that in view of Lemma \ref{ext lemma} and Lemma \ref{tor lemma}, we need to study the behaviour of the local $v$-number of the graded $R$-module $(U+I^nV)/I^nW$, i.e., of the graded $R$-module $M_n'$. Note that with the hypothesis above and by \cite[Proposition 3.4]{KW04} and its proof, the set $\Ass_R \left(U+I^nV/I^nW\right)$ stabilizes eventually, that is the set $\Ass_R(M_n')$ stabilize for all $n\gg0.$

\section{Associated primes of \texorpdfstring{$(U+I^nV)/I^nW$}{(U+InV)/InW}}

In this section, under the Hypothesis \ref{hypo}, we investigate the behavior of the set $\Ass_R \left(U+I^nV/I^nW\right)$ and establish a result analogous to \cite[Theorem 2.2]{GNP26}.
%Before that we need to fix some notations.
% \begin{notation}\label{notation}
% For each integer $n\geq 1$, put
% $M_n=U+I^nV, D_n=
% M_n\cap\left((0:_E I)+I^{n-1}W\right), M_n'=M_n/I^nW, D_n'=D_n/I^nW,$ and $T_n=(M_n\cap I^{n-1}W)/I^nW$. Note that $I^nW\subseteq M_n$ as $W\subseteq V$. Also since $I^nW\subseteq D_n$, 
% $D_n'$ is well-defined. 
% \end{notation}
% We write $V(I)$ for the set of prime ideals
% of $R$ containing $I$. We begin with the following lemma.

\begin{lemma}\label{lem1}
Assume the Hypothesis as in \ref{hypo}.
There exists $k\in\mathbb{N}$ such that for every
$n\geq k$, the following statements hold:
\begin{enumerate}[label=(\alph*), font=\upshape]
    \item For every $\p\in V(I)$, we have
\[ \left\{ x\in M_n : \p=(I^nW:_R x) \right\}
\subseteq
D_n.\]

\item $\Ass_R(M_n')\cap V(I)=\Ass_R(D_n').
$

\end{enumerate}
\end{lemma}

    \begin{proof}
Since $W\subseteq E$, by \cite[Corollary  2.1]{Dao21} there exists
$k\in\mathbb{N}$ such that for all $n\geq k$
\begin{equation}\label{(1)}
(I^nW:_E I) = (0:_E I)+I^{n-1}W
\end{equation}

\textnormal{(a)}
Let $n\geq k$ and let $\p\in V(I)$.
Suppose that $x\in M_n$ satisfies
$\p=(I^nW:_R x).$
Since $\p\in V(I)$, we have
$I\subseteq\p$. Consequently,
$Ix\subseteq I^nW.$
Thus,
$x\in(I^nW:_E I).$
 By (\ref{(1)}) we obtain
$x\in(0:_E I)+I^{n-1}W,$ for all $n\geq k$ .
Since $x\in M_n$, it follows that
$x\in
M_n\cap\left((0:_E I)+I^{n-1}W\right)=D_n$, for all $n\geq k$
Therefore,
$$\left\{
x\in M_n:
\p=(I^nW:_R x)
\right\}
\subseteq
D_n$$
for every $n\geq k$.
\medskip

\textnormal{(b)}
For every $n\geq 1$, 
$D_n=
M_n\cap\left((0:_E I)+I^{n-1}W\right)\subseteq M_n$. Therefore,
$
D_n'
\subseteq
M_n'$, and hence 
\begin{equation}\label{(2)}
\Ass_R
(D_n')
\subseteq
\Ass_R
(M_n').
\end{equation}

Note that $ID_n\subseteq IM_n\cap \left(I(0:_EI)+I^nW\right)=IM_n\cap I^nW\subseteq I^nW$. Consequently,
$
ID_n'=I\left({D_n}/{I^nW}\right)=\left(ID_n+I^nW\right)/I^nW=0.$ Hence,
\begin{equation}\label{(3)}
  \Ass_R
(D_n')
\subseteq V(I).
 \end{equation}
 Combining equation (\ref{(2)}) and (\ref{(3)}) we obtain
 \begin{equation}\label{(4)}
  \Ass_R(D_n')
\subseteq V(I)\cap \Ass_R(M_n').
 \end{equation}
Conversely, let $n\geq k$ and $\p\in V(I)\cap \Ass
(M_n')$. By definition of associated primes, there exists an element $\overline{x}\in M_n'={M_n}/{I^nW}$ such that 
$\p= (0:_R\overline x)$ = $(I^nW:_R x)$. 

Since $\p\supseteq I$, it follows
from part \textnormal{(a)} that
$x\in
M_n\cap\left((0:_E I)+I^{n-1}W\right)
=D_n.$
Consequently,
$\overline{x}\in{D_n}/{I^nW}=D_n'.$
Thus,
$\p\in
\Ass_R
(D_n').
$
Therefore,
\begin{equation}\label{(5)}
\Ass_R
(M_n')
\cap V(I)
\subseteq
\Ass_R
(D_n').
\end{equation}

Combining
(\ref{(4)}) and
(\ref{(5)}), we obtain part (b).
\end{proof}

We now prove the main result of this section.
\begin{theorem}\label{ass prime th Mn' and Tn}
    Assume the Hypothesis as in \ref{hypo}. we have, for all sufficiently large integers $n$,
\[
\Ass_R(M_n')\cap V(I)
=
\Ass_R(0:_UI)\cup\Ass_R(T_n).
\]
\end{theorem}

\begin{proof}
    Let $C=(0:_E I)=\{x\in E:Ix=0\}.$  Applying the Artin-Rees lemma  to the inclusion $C\subseteq E$, there exists an integer $r$ such that $C\cap I^nE = I^{n-r}(C\cap I^rE)$ for all $n\geq r+1$.
Since $IC=0$, it follows that
$C\cap I^nE=0$
for all sufficiently large $n$.
In particular,
\begin{equation}\label{(6)}
C\cap I^{n-1}W=0
\end{equation}
for all sufficiently large $n$. We claim that
\begin{equation}\label{(7)}
D_n=(C\cap M_n)\oplus
(M_n\cap I^{n-1}W)
\end{equation}
for all sufficiently large $n$.

Indeed, let $x\in D_n$. By the definition of
$D_n$, we have $x\in M_n$ and $x\in C+I^{n-1}W.$
Thus, we can write $x=a+b,$ where $a\in C$ and $b\in I^{n-1}W$.
We first claim that $a\in M_n$. Let "$\overline{(-)}$" denote reduction modulo $U$. Consider the quotient module $\overline{E}=E/U,$ and let $\overline{C}$ denote the image of $C$ in $\overline{E}$. Since $IC=0$, we have $I\overline{C}=0.$ Applying the Artin-Rees lemma   to the inclusion $\overline{C}\subseteq\overline{E}$, we obtain
$\overline{C}\cap I^n\overline{E}=0$ for all sufficiently large $n$.
Now,
$\overline{x}=\overline{a}+\overline{b}$ and hence
$\overline{a}=\overline{x}-\overline{b}.$ Since
$x\in M_n=U+I^nV\subseteq U+I^nE,$
we have $\overline{x}\in I^n\overline{E} \subseteq I^{n-1}\overline{E}.$ Also, as $b\in I^{n-1}W$, we have
$\overline{b}\in I^{n-1}\overline{E}.$
Consequently, $\overline{a}\in I^{n-1}\overline{E}.$
On the other hand, $a\in C$ implies that $\overline{a}\in\overline{C}$. Therefore, $\overline{a}\in
\overline{C}\cap I^{n-1}\overline{E}=0$ for all sufficiently large $n$. Thus $\overline{a}=0$, which implies that
$a\in U\subseteq M_n.$
Since $x,a\in M_n$, it follows that
$b=x-a\in M_n.$ Therefore,
$b\in M_n\cap I^{n-1}W.$ Hence,
$x\in (C\cap M_n)+(M_n\cap I^{n-1}W),$ and hence,
$D_n\subseteq (C\cap M_n)+(M_n\cap I^{n-1}W).$

Conversely, by the definition of $D_n$, we have
$C\cap M_n\subseteq D_n$ and $M_n\cap I^{n-1}W\subseteq D_n.$
Thus, $D_n=(C\cap M_n)+(M_n\cap I^{n-1}W).$
Moreover, by (\ref{(6)}), $(C\cap M_n)\cap(M_n\cap I^{n-1}W) \subseteq C\cap I^{n-1}W=0.$ Therefore, $D_n=(C\cap M_n)\oplus (M_n\cap I^{n-1}W).$

Now define an $R$-module homomorphism
$$\varphi:(C\cap M_n)\oplus T_n
\longrightarrow {D_n}/{I^nW} \qquad \text{by} \qquad
\varphi(a,\overline{b})=\overline{a+b},$$
where $a\in C\cap M_n$ and
$b\in M_n\cap I^{n-1}W$.

It is clear that, $\varphi$ is a well-defined
$R$-module homomorphism. To prove injectivity, suppose that
$\varphi(a,\overline{b})=0.$
Then $a+b\in I^nW \subseteq I^{n-1}W.$
Since $b\in I^{n-1}W$, we obtain
$a\in C\cap I^{n-1}W=0,$ for all sufficiently large $n$.
Thus $a=0$, and consequently,
$b\in I^nW$. Therefore,
$(a,\overline{b})=(0,\overline{0}),$ which proves that $\varphi$ is injective for $n\gg 0$. To prove surjectivity, let
$\overline{x}\in D_n/I^nW$. We can write
$x=a+b,$ where
$a\in C\cap M_n
\text{ and }
b\in M_n\cap I^{n-1}W.$
Hence
$\varphi(a,\overline{b}) =\overline{a+b} =\overline{x}.$
Thus $\varphi$ is surjective.
Consequently, $${D_n}/{I^nW} \simeq (C\cap M_n)\oplus T_n, \qquad  n\gg 0 .$$
By the definitions of $D_n'$ and $T_n$,
we obtain
$D_n'\simeq (C\cap M_n)\oplus T_n$, for all sufficiently large $n$.
Therefore,
$\Ass_R(D_n') = \Ass_R(C\cap M_n)\cup\Ass_R(T_n).$
Since $C\cap M_n=(0:_{M_n}I),$ it follows that
$\Ass_R(D_n') = \Ass_R((0:_{M_n}I))\cup\Ass_R(T_n).$
Finally, by Lemma \ref{lem1}(b), $\Ass_R(M_n')\cap V(I)=\Ass_R(D_n')$ for all sufficiently large $n$.
Hence
$\Ass_R(M_n')\cap V(I)
=
\Ass_R(0:_{M_n}I)\cup\Ass_R(T_n).
$
This completes the proof using the fact $(0:_{M_n}I)=(0:_UI)$ for $n\gg 0$, see \cite[Lemma 3.7]{GP25}.
\end{proof}

We end this section with the following useful remark.

\begin{remark}\label{rmk1}
Note that in Theorem \ref{ass prime th Mn' and Tn} we have proved for $n\gg 0,$
$$D_n'\simeq (C\cap M_n)\oplus T_n\simeq (0:_{U+I^nV}I)\oplus T_n\simeq (0:_{U}I)\oplus T_n$$
The last isomorphism follows from \cite[Lemma 3.7]{GP25}. Let $\overline{(-)}$ denotes modulo $\Gamma_I(U)$. Then $\overline{D_n'}\simeq\overline{T_n}$ since $(0:_{\overline{U}}I)=0$. Also, by definition we have 

\begin{align*}
D_n'
&= {D_n}/{I^nW} 
\\ &= {M_n \cap \big((0 :_E I) + I^{n-1}W\big)}/{I^nW} 
\\ &= {M_n \cap (I^nW :_E I)}/{I^nW}
\qquad \text{ for $n\gg 0$, by \cite[Corollary 2.1]{Dao21}}
\\ &= {(I^nW :_{M_n} I)}/{I^nW}
\qquad \text{ for $n\gg 0$},
\\ &= (0 :_{M_n/I^nW} I)
\qquad \text{ for $n\gg 0$},
\\ &= (0 :_{M_n'} I) 
\qquad \text{ for $n\gg 0$}.
\end{align*}
It follow that $(0 :_{\overline{M_n'}} I)\simeq\overline{T_n} \, $\text{ for $n\gg 0$} .
\end{remark}
\section{V-numbers of Ext and Tor modules}
Let the hypothesis be as in \ref{hypo}. In this section, we first prove that the local $v$-number of $M'_n$ is eventually a non-zero constant or a linear polynomial. We also prove that if $(0:_UI)\neq 0$ then the global $v$-number of $M'_n$ is eventually equal to a non-zero constant provided  $I$ is generated in positive degree and  $\Ass_R (M'_n) \subseteq V(I)$. Finally, using the identification of $\Ext^k_R(L,I^nM/N)$ and  $\Tor_k^R(L,I^nM/N)$ (see Lemma \ref{ext lemma} and \ref{tor lemma}), we prove the main result (Theorem \ref{main result}) of this article. 
 We begin with the following preparatory lemma.
\begin{lemma}\label{lem2}
Assume the Hypothesis as in \ref{hypo}. Let $\overline{(-)}$ denote reduction modulo $\Gamma_I(U)$.
Then $T_n\simeq\overline{T_n}$ for $n\gg 0$.
\end{lemma}

\begin{proof}
Note that
$T_n={(M_n\cap I^{n-1}W})/{I^nW}.$
First, we show that $\Gamma_I(U)\cap I^{n-1}W=0$
for $n\gg 0$. Since $U$ is finitely generated,
there exists an integer $k\geq 1$ such that
$I^k\Gamma_I(U)=0$. Applying  Artin-Rees lemma on $\Gamma_I(U)\cap W \subset W$,
there exists an integer $c\geq 0$ such that, for all sufficiently large $n$,
$$
\Gamma_I(U)\cap I^{n-1}W
=I^{n-1-c}
  \bigl(\Gamma_I(U)\cap I^cW\bigr)\\
\subseteq I^{n-1-c}\Gamma_I(U)\\
=0.$$
Therefore, for all $n\gg 0$,
\begin{align*}
\overline{T_n}
&={(M_n\cap
I^{n-1}W)+(\Gamma_I(U)+I^nW)}/{I^nW+\Gamma_I(U)}\\
&\simeq M_n\cap I^{n-1}W/(M_n\cap
I^{n-1}W)\cap ({I^nW+\Gamma_I(U)})\\
&=
 M_n\cap I^{n-1}W/I^nW+
{(\Gamma_I(U)\cap I^{n-1}W)} \\
&={M_n\cap I^{n-1}W}/
{I^nW}\\
&=T_n.
\end{align*}
Here, the second equality follows from the second isomorphism theorem and the third equality follows from the modular law since $I^nW\subseteq M_n\cap
I^{n-1}W$. The second last equality holds from the fact that $\Gamma_I(U)\cap I^{n-1}W=0$ for $n\gg 0$.
\end{proof}
 Note that by Theorem \ref{ass prime th Mn' and Tn}, we have for all $n\gg0,$ $\Ass_R (M_n')\cap V(I)=\Ass_R (0:_U I)\cup\Ass_R (T_n)$. In the next result, we consider the case when $\p\in \Ass (M_n'\cap V(I))\setminus \Ass (0:_UI)$ and prove the asymptotic linearity of $v_\p(M_n')$ for all $n\gg0.$

\begin{theorem}\label{if p is not ass prime of 0 collon I}
    Assume the Hypothesis as in \ref{hypo}. Let $\p\in  \left(\Ass_R(M_n')\cap V(I)\right)\setminus \Ass (0:_UI)$. Then there exists $a\in\{d_1,\ldots,d_c\}$ and $b\in \Z$ such that $v_{\p}(M_n')=a n+b$ for $n\gg 0$.
\end{theorem}
\begin{proof}
    Note by Theorem \ref{ass prime th Mn' and Tn} we have for all $n\gg0,$ $$\Ass_R (M_n')\cap V(I)=\Ass_R (0:_U I)\cup\Ass_R (T_n).$$ 
    By assumption, $\p \in \Ass_R (T_n)$. Note that $\overline{M_n'}=(U+I^nV)/(I^nW+\Gamma_I(U))\simeq \overline{U}+I^n{\overline{V}}/I^n{\overline{W}}$, where $\overline{(-)}$ denotes modulo $\Gamma_I(U)$. Consider the following short exact sequence 
    \begin{equation*}
        0\to I^nW+\Gamma_I(U)/I^nW\to U+I^nV/I^nW\to U+I^nV/I^nW+\Gamma_I(U)\to 0.
    \end{equation*}
    By Artin-Rees lemma, $I^n W\cap \Gamma_I(U)=0$ for $n\gg 0$. Consequently, we get $I^nW+\Gamma_I(U)/I^nW=\Gamma_I(U)/\Gamma_I(U)\cap I^nW=\Gamma_I(U)$ for $n\gg 0$. Therefore, for all $n\gg0$, the above exact sequence becomes 
    \begin{equation}\label{equation of Mn and Mn'}
        0\to \Gamma_I(U)\to M_n'\to \overline{M_n'}\to 0.
    \end{equation}
    This gives us $\Ass_R (M_n')\subseteq \Ass_R(\Gamma_I(U))\cup \Ass_R (\overline{M_n'}),$  for all sufficiently large $n$. Lemma \ref{ass prime of gammaI} implies that $\p\in \Ass (\overline{M_n'})$. Since $(0:_{\overline{U}}I)=0$, by \cite[Theorem 3.8]{GP25}  there exists $a\in\{d_1,\ldots,d_c\}$ and $b\in \Z$  such that
    $$v_{\p} (\overline{M_n'})=a n+b   \text{ for}\ n\gg 0.$$

    We claim that $v_{\p} (M_n')=v_{\p} (\overline{M_n'}).$ By Proposition \ref{v number in exact seq}, $v_{\p} (\overline{M_n'})\leq v_{\p} (M_n')\leq v_{\p} (T_n)  .$
    Also by Lemma \ref{lem2}, we have $v_{\p} (T_n)=v_{\p} (\overline{ T_n})$. By Remark \ref{rmk1}, $(0:_{\overline{M'_n}} I)\simeq \overline{T_n} \quad \text{for $n\gg 0$}.$ Let
$\p = (0:_R \overline{f}) $ for some $\overline{f}\in\overline{M'_n},$ such that $
\deg(\overline{f})=v_{p}(\overline{M'_n}).$ Since $\p \in V(I)$, we have $I\subseteq{\p}=(0:_R\overline{f}).$ Thus, $I\overline{f}=0,$ and hence
$\overline{f}\in(0:_{\overline{M_n'}}I)=\overline{T_n}.$
Therefore, $v_{\p}(\overline{T_n})\leq \deg(\overline{f})
=v_{\p}(\overline{M_n'})$ for all $n\gg 0$. This implies that $v_{\p}(M_n')=v_{\p}(\overline{M_n'})=an+b$ for all $n\gg0.$ 
\end{proof}

In the next result, we deal with the case when $\p\in \Ass_R ((0:_U I))$. Since the set $\Ass_R (0:_U I)\cap\Ass_R (T_n)$ may be non-empty, we split the proof in two cases depending on whether $\p\in \Ass_R(T_n)$ or $\p\notin \Ass_R(T_n)$. We begin with the case when $\p\notin \Ass_R(T_n)$.
\begin{theorem}\label{if p is not ass prime of Tn}
   Assume the Hypothesis as in \ref{hypo}. Let $\p\in  (\Ass_RM_n'\cap V(I))\setminus \Ass (T_n)$. Then, $v_{\p}(M_n')=v_{\p}(0:_UI)$ for $n\gg 0$.
\end{theorem}    
\begin{proof}
By Theorem \ref{ass prime th Mn' and Tn}, for all sufficiently large
integers $n$, we have
\[
\Ass_R(M_n')\cap V(I)
=
\Ass_R(0:_UI)\cup\Ass_R(T_n).
\]    
Since $\p\notin \Ass_R (T_n)$, $\p\in \Ass_R((0:_UI))$. By Remark \ref{rmk1}, we have $(0 :_{M_n'} I)=(0:_UI)\oplus T_n$ for all $n\gg0$. This gives the following two exact sequences
\begin{equation}\label{1st seq in 2nd main}
 0\to (0:_UI)\to (0 :_{M_n'} I)\to T_n\to 0   
\end{equation}
 and 
 \begin{equation}\label{2nd seq in 2nd main}
     0\to T_n\to (0 :_{M_n'} I)\to (0:_UI)\to 0.
 \end{equation}
 Using Proposition \ref{v number in exact seq} for the exact sequences (\ref{1st seq in 2nd main}) and (\ref{2nd seq in 2nd main}), we obtain $v_{\p}(0 :_{M_n'} I)=v_{\p}(0 :_U I).$ The result now follows from Proposition \ref{v number of M and colon same}.
\end{proof}
Now we address the case when $\p \in \Ass_R(T_n)$ in the following theorem. However, to establish a comparable result, we must assume that $I$ is generated in positive degrees.
\begin{theorem}\label{if p is ass prime of Tn}
     Assume the Hypothesis as in \ref{hypo}. Let $\p\in \Ass_R(0:_UI)\cap\Ass_R(T_n)$. Suppose that $d_i\geq 1$ for $1\leq i\leq c$. Then,  $v_{\p}(M_n')=v_{\p}(0:_UI)$ for $n\gg 0$.
\end{theorem}
\begin{proof}
Due to Remark \ref{rmk1}, we have $(0 :_{M_n'} I)\simeq(0:_UI)\oplus T_n$ for all $n\gg0$. By Lemma \ref{easy} we have, $v_{\p}(M_n')=\operatorname{min}\{v_{\p}(0:_UI), v_{\p}(T_n)\}.$ Since $T_n=(M_n\cap I^{n-1}W)/I^nW$ is submodule of $I^{n-1}W/I^nW$, by Proposition \ref{v number in exact seq}, $v_{\p}(I^{n-1}W/I^nW)\leq v_{\p}(T_n)$. By \cite[Theorem 2.11]{FG25}, $v_{\p}(I^{n-1}W/I^nW)=an+b$ for $n\gg 0$, where $a\in\{d_1,\ldots,d_c\}$ and $b\in \Z$. Since $a\geq 1$, we have $v_{\p}(T_n) \to \infty $ as $n\to \infty$. Thus, for all $n\gg 0$, $v_{\p}(T_n)>v_{\p}(0:_UI)$. This proves the result.
\end{proof}

\begin{theorem}\label{whn p is in intersection of prime and colon}
    Assume the Hypothesis as in \ref{hypo}. Let $\p\in \Ass_R (M_n')\cap \Ass_R(0:_UI)$. If $d_i\geq 1$ for $1\leq i\leq c$, then $v_{\p}(M_n')=v_{\p}(0:_UI)$ for $n\gg 0$.
\end{theorem}
\begin{proof}
By Theorem \ref{ass prime th Mn' and Tn}, for all sufficiently large
integers $n$, we have
\[
\Ass_R(M_n')\cap V(I)
=
\Ass_R(0:_UI)\cup\Ass_R(T_n).
\]     
There are two possibilities: either $\p\in \Ass_R (T_n)$ or $\p\notin \Ass_R (T_n)$. The result now follows from Theorem \ref{if p is not ass prime of Tn} and Theorem \ref{if p is ass prime of Tn}.
\end{proof}
Next, we describe the asymptotic behavior of the global Vasconcelos invariant of $M_n'$ when the module $(0:_UI)$ is non-zero.

\begin{theorem}\label{result on global v number}
Assume the Hypothesis as in \ref{hypo}. Let $(0:_UI)\neq 0$. Then, $v(M_n')\leq v(0:_UI)$ for $n\gg 0$. Further if $d_i\geq 1$ for $1\leq i\leq c$ and $\Ass_R (M_n')\subseteq V(I)$, then for all $n\gg0$, $$v(M_n')= v(0:_UI).$$
\end{theorem}
\begin{proof}
Due to Remark \ref{rmk1}, we have $(0 :_{M_n'} I)=(0:_UI)\oplus T_n$ for all $n\gg0$. Therefore, $(0:_UI)$ is isomorphic to a submodule of $M_n'$. Thus, by Proposition \ref{v number in exact seq}, $v_{\p}(M_n')\leq v_{\p}(0:_UI)$ whenever $\p\in \Ass_R ((0:_UI))$. 
Note that \begin{align*}
  v(M_n')&=\inf\{v_{\p}(M_n')\mid \p\in \Ass_R (M_n')\}\\ &\leq \inf\{v_{\p}(M_n')\mid \p\in \Ass_R (0:_UI)\}\\ &\leq \inf\{v_{\p}(0:_UI)\mid \p\in \Ass_R (0:_UI)\}\\&= v(0:_UI).
\end{align*} 
Now by assumption $\Ass_R (M_n')\subseteq V(I)$. Let $\p\in \Ass_R (M_n')$. Note that by Theorem \ref{ass prime th Mn' and Tn} we have $$\Ass_R (M_n')\cap V(I)=\Ass_R (0:_U I)\cup\Ass_R (T_n).$$ 
This implies that $\p\in \Ass_R (0:_U I)\cup\Ass_R (T_n).$ If $\p\notin \Ass_R (0:_U I)$, then by Theorem \ref{if p is not ass prime of 0 collon I}, $v_{\p}(M_n')=an+b\to \infty$ as $n\to \infty$. In other cases, by Theorem \ref{if p is not ass prime of Tn} and Theorem \ref{if p is ass prime of Tn}, we have $$v_{\p}(M_n')= v_{\p}(0:_UI).$$ This proves that $v(M_n')= v(0:_UI).$
\end{proof}
We conclude this section by combining the preceding results to establish the main theorem of this article, stated as Theorem \ref{main result} in the introduction.
\begin{theorem}\label{proof of main result}
   Let $R$ be a Noetherian $\N$-graded ring. Assume that $L$ and  $M$ are finitely generated $\Z$-graded $R$-modules, and $I$ a homogeneous ideal of $R$,  generated by homogeneous elements $y_1,\ldots,y_c$ of degree $d_1\leq\ldots\leq d_c$ respectively. Let $N$ be a graded submodule of $M$. Fix $k\in \N$. Denote $H:=\Ext^k_R(L,M)$ (or $\Tor_k^R(L,M)$) and $H_n:=\Ext^k_R(L, M/I^nN)$ (or $\Tor_k^R(L,M/I^nN)$). Then the following statements hold:
   \begin{enumerate}[label=(\arabic*), font=\upshape]
       \item  Let $\p\in \left(\Ass_R H_n\cap V(I)\right)\setminus \Ass_R(0:_HI)$ for $n\gg 0$. Then there exists $a_\p\in\{d_1,\ldots,d_c\}$ and $b_\p\in \Z$ such that $v_{\p}(H_n)=a_\p n+b_\p$ for $n\gg 0$.
       \item  Let $\p\in \Ass_R(H_n)\cap \Ass_R(0:_HI)$ for $n\gg 0$. If $d_i\geq 1$ for $1\leq i\leq c$, then $v_{\p}(H_n)=v_{\p}(0:_HI)$ for $n\gg 0$.
       \item Let $(0:_HI)\neq 0$. Then, $v(H_n)\leq v(0:_HI)$ for $n\gg 0$. Further, if $d_i\geq 1$ for $1\leq i\leq c$, and $\Ass_R(H_n)\subseteq V(I))$, then $v(H_n)= v(0:_HI)$ for $n\gg 0$. 
   \end{enumerate}
\end{theorem}
\begin{proof}
We first prove (1) for the $\Ext$-module and the similar proof will follow for the $\Tor$-modules as well. Note that by Lemma \ref{ext lemma}, there exists $q\in \N$ such that
    $$\Ext^k_R(L,M/I^nN)\cong (U+I^{n-q}V)/I^{n-q}W=M'_{n-q}  \ \ \text{for all}\ n\geq q$$
    and for some graded submodules $U,V$ and $W$ of a finitely generated $R$-module $E$ satisfying $W\subseteq V$ and $U\cong \Ext^k_R(L,M)$. From the given condition, one has $\p\in \Ass (M_n'\cap V(I))\setminus \Ass ((0:_UI))$ for $n\gg 0$. By Theorem \ref{if p is not ass prime of 0 collon I}, there exists $a_\p\in\{d_1,\ldots,d_c\}$ and $b\in \Z$ such that $ v_\p \left((U+I^{n}V)/I^{n}W\right)=a_\p n+b$ for $n\gg 0$. Set $b_p=b-a_\p q$. Therefore,
    $$v_\p\left(\Ext^k_R(L,M/I^nN)\right)= v_\p \left((U+I^{n-q}V)/I^{n-q}W\right)=a_\p n+b_\p \ \text{for all}\ n\gg 0.$$

Similarly, (2) and (3) follow from Theorem \ref{whn p is in intersection of prime and colon} and Theorem \ref{result on global v number}, respectively.
\end{proof}
\section{Examples}
In this section, we give some examples in support of our results. We start
with the following which can be derived from \cite[Example 3.1]{FG25} as observed in \cite[Example 4.1]{GP25}.
\begin{example}
    Let $R=K[X,Y]$ be a standard graded polynomial ring in two variables $X$ and $Y$ over a field $K$. Consider $U=V=W:=R/(XY^b)$, $I=(X^{a})$, $\p:=(X)$, and $\mathfrak{m}:=(X,Y)$,
where $a$ and $b$ are some positive integers. Then
\begin{enumerate}
    \item $A=(0:_U I)=Y^bU\neq 0$. Thus, $\Ass_R(A)=\{\p\}$ and $\p=(XY^b:_RY^b).$ Therefore, $v_{\p}(A)=b.$
    \item $\Ass_R(U+I^nV/I^nW)=\{\p,\m\}$ if $bn\geq 2$ and $\p$ if $bn=1.$
    \item $v_{\p}(U+I^nV/I^nW)=b$ whenever $an\geq 2.$
    \item $v_{\m}(U+I^nV/I^nW)=an+(b-2)$ whenever $bn\geq 2.$
    \item $v(U+I^nV/I^nW)=b$ whenever $an\geq 2.$
\end{enumerate}
In particular, whenever $(0:_UI)\neq 0$ and $\p\in \Ass_R(0:_UI)$, $v_{\p}(U+I^nV/I^nW)=v(U+I^nV/I^nW)=v_{\p}(0:_UI)$ for all $n\gg0$. If $\p\in \Ass_R(U+I^nV/I^nW)\setminus \Ass_R(0:_UI)$ for all $n\gg0$, then $v_{\p}(U+I^nV/I^nW)$ is eventually linear.
\end{example}

Recall that, by Proposition \ref{easy}, $v(M(r))=v(M)-r$ and $v_{\p}(M(r))=v_{\p}(M)-r$ for all $\p \in \Ass_R(M)$ and $r\in \Z.$ We use this result repeatedly in the next example, which is the general version of \cite[Example 4.4]{GP25} and satisfies Theorem \ref{proof of main result}.

\begin{example}
    Let $R = K[X,Y]/(X^aY^b)$ with $a,b\geq 1$ and $\deg(X)=\deg(Y)=1,$ where $K$ is a field. Suppose $x$ and $y$ are the residue classes of $X$ and $Y$ in $R$, respectively. Set $L= R/(x^a)$, $I=M=N=(y^b)$, and $\m=(x,y)$. Then,
    \begin{enumerate}
        \item $\displaystyle (0:_{\Ext^k_R(L,M)}I)=\begin{cases}
            0, & \text{if } k =0 \text{ or } k\geq 1 \text{ is odd }\\ (y^b)/(y^{2b})(r(a+b)), & \text{if } k=2r\geq 2.
        \end{cases}$\\
        For $k=2r\geq 2$, $\displaystyle\Ass_R\left((0:_{\Ext^k_R(L,M)}I)\right)=\{\m\}$ and $\displaystyle v\left((0:_{\Ext^k_R(L,M)}I)\right)=(1-r)(a+b)+b-2.$

       \item $\displaystyle (0:_{\Tor_k^R(L,M)}I)=\begin{cases}
            0, & \text{if } k \geq 0 \text{ even}\\ (y^b)/(y^{2b})(-((r+1)a+rb)), & \text{if } k=2r+1\geq 1.
        \end{cases}$\\
For $\displaystyle k=2r+1\geq 1$, $\displaystyle\Ass_R\left((0:_{\Tor_k^R(L,M)}I)\right)=\{\m\}$ and $\displaystyle v\left((0:_{\Tor_k^R(L,M)}I\right)=(r+2)(a+b)-2$.

       \item $\displaystyle\Ass_R\left(\Ext^k_R(L,M/I^nN)\right)=\Ass_R\left(\Tor^k_R(L,M/I^nN)\right)=\{\m\}$ for all $k\geq 0$ and $n\geq 1.$
       
        \item For all $n\geq 1$,
        
        $\displaystyle v\left(\Ext^k_R(L,M/I^nN)\right)=\begin{cases}
            nb+a+b-2 , & \text{if } k =0 \\ nb-r(a+b)+b-2 & \text{if } k=2r+1\geq 1 \\ (1-r)(a+b)+b-2, & \text{if } k=2r\geq 2.
        \end{cases}$

        \item For all $n\geq 1$,
        
        $\displaystyle v\left(\Tor_k^R(L,M/I^nN)\right)=\begin{cases}
            nb+a+b-2 , & \text{if } k =0 \\ (r+2)(a+b)-2 & \text{if } k=2r+1\geq 1 \\ nb+ (r+1)(a+b)-2, & \text{if } k=2r\geq 2.
        \end{cases}$
    \end{enumerate}
\end{example}

\begin{proof}
    Note that $$\mathbb F_L: \ldots \xrightarrow{y^b}R(-(2a+b))\xrightarrow{x^a}R(-(a+b))\xrightarrow{y^b}R(-a)\xrightarrow{x^a}R\to 0$$ is the minimal graded free resolution of $L.$ Applying the functors $\Hom_R(-,M)$ and $\otimes_RM$ to the above resolution, we obtain the following complexes
    \begin{equation}\label{complex}
        \begin{split}
            \Hom_R(\mathbb F_L,M) : &\, 0\to M\xrightarrow{x^a}M(a)\xrightarrow{y^b}M(a+b)\xrightarrow{x^a}M(2a+b)\xrightarrow{y^b} \ldots \\ \mathbb F_L\otimes_R M : &\, \ldots \xrightarrow{y^b} M(-(2a+b))\xrightarrow{x^a} M(-(a+b)) \xrightarrow{y^b}M(-a)\xrightarrow{x^a}M\to 0.
        \end{split}
    \end{equation}
The shifts of the modules in the above complexes follow from the fact that $\Hom_R(R(-k),W)\cong W(k)$ and $R(-k)\otimes_R W\cong W(-k)$ for any graded $R$-module $W.$

Since $M=(y^b)$, from \eqref{complex}, we get 

\begin{equation}\label{e}
    \Ext^k_R(L,M)=\begin{cases}
            M, & \text{if } k=0 \\ 0, & \text{if } k \text{ odd}\\ (y^b)/(y^{2b})(r(a+b))  & \text{if } k=2r\geq 2.
         \end{cases}
\end{equation}

\begin{equation}\label{t}
 \text{and }    \Tor_k^R(L,M)=\begin{cases}
              M, & \text{if } k=0 \\ (y^b)/(y^{2b})(-((r+1)a+rb)), & \text{if } k=2r+1\geq 1 \\ 0,  & \text{if } k \text{ even}.
         \end{cases}
\end{equation}
Since $(x^a,y^b)$ annihilates modules in \eqref{e} for all $k=2r\geq 2$, we have  $\Ass_R((0:_{\Ext^k_R(L,M)}I))=\{\m\}$. Further, $\m=(\overline{0}:_R\overline{z})$, where $\overline{z}=\overline{y^{2b-1}x^{a-1}}\in (y^b)/(y^{2b}).$ Therefore, $v((0:_{\Ext^k_R(L,M)}I)=a+2b-2-r(a+b)=b-2-(r-1)(a+b).$ This proves (1). From \eqref{t} and by the same argument as above, (2) follows.

Since $M/I^nN= (y^b)/(y^{(n+1)b})$ for all $n\geq 1$ and by replacing $M$ by $M/I^nN$ in \eqref{complex}, we obtain

\begin{equation}\label{1}
    \Ext^k_R(L,M/I^nN)=\begin{cases}
            (y^b)/(y^{(n+1)b}), & \text{if } k=0 \\ (y^{nb})/(y^{(n+1)b})((r+1)a+rb), & \text{if } k=2r+1\geq 1\\ (y^b)/(y^{2b})(r(a+b))  & \text{if } k=2r\geq 2.
         \end{cases}
\end{equation} 
\begin{equation}\label{2}
   \text{and }  \Tor_k^R(L,M/I^nN)=\begin{cases}
              (y^b)/(y^{(n+1)b}), & \text{if } k=0 \\ (y^b)/(y^{2b})(-((r+1)a+rb)), & \text{if } k=2r+1\geq 1 \\ (y^{nb})/(y^{(n+1)b})(-(r(a+b))),  & \text{if } k=2r\geq 2.
         \end{cases}
\end{equation}
 Note that $(x^a,y^{nb})$ or $(x^a,y^b)$ annihilates the modules in \eqref{1} and \eqref{2}. Therefore, we obtain  $$\Ass_R\left(\Ext^k_R(L,M/I^nN)\right)=\Ass_R\left(\Tor_k^R(L,M/I^nN)\right)=\{\m\}$$ for all $k\geq 0$ and $n\geq 1$. This proves (3).

We now prove (4). By \eqref{1}, we have for $k=0 \text{ and }2r+1\geq 1$, we have $\m=(\overline{0}:_R \overline{z})$, where $\overline{z}=\overline{y^{(n+1)b-1}x^{a-1}}\in (y^b)/(y^{(n+1)b})$. Therefore, for $k=0,$ $v\left(\Ext^k_R(L,M/I^nN)\right)=nb+a+b-2$. For $k=2r+1\geq 1$, we also have $\m=(\overline{0}:_R \overline{z})$, where $\overline{z}=\overline{y^{(n+1)b-1}x^{a-1}}\in (y^{nb})/(y^{(n+1)b})$. Thus, $v(\Ext^k_R(L,M/I^nN))=(n+1)b+a-2-(r+1)a-rb=nb-r(a+b)+b-2$ for all $n\geq 1.$ Now, for $k=2r\geq 2$, we have $\m=(\overline{0}:_R\overline{z})$, where $\overline{z}=\overline{y^{2b-1}x^{a-1}}\in (y^{b})/(y^{2b})$. Therefore, $v( \Ext^k_R(L,M/I^nN))=2b+a-2-(r(a+b))=(1-r)(a+b)+b-2=v((0:_{\Ext^k_R(L,M)}I))$ for all $n\geq 1.$

Using \eqref{2} and by the same argument as above, (5) follows.
\end{proof} 
\noindent
\section*{Acknowledgment}
 This work was initiated when the three authors met at the Advanced Instructional School (AIS) on Topics in Commutative Algebra and Algebraic Geometry, held at IIT Palakkad in June 2026. The authors gratefully acknowledge the institute for its warm hospitality and wonderful campus environment. The first named author acknowledges the Government of India for support through the Prime  Minister's Research Fellowship (ID: 1303161). The second named author acknowledges the University Grants Commission (UGC) for their support (Reference number 231610213879).

\bibliographystyle{alpha}
 
\end{document}